\documentclass[12]{article}

\usepackage{amssymb}
\usepackage{color, amsmath,amssymb, amsfonts, amstext,amsthm, latexsym}

\usepackage{graphicx}
\usepackage{longtable}

\renewcommand{\phi}{\varphi}

\newcommand{\la}{{\lambda}}
\newtheorem{theorem}{Theorem}

\newtheorem{lemma}{Lemma}
\newtheorem{definition}{Definition}

\def\EE{{\mathbb{E}}}

\begin{document}

\title{Higher-order Derivative Local Time for
Fractional Ornstein-Uhlenbeck Processes
\thanks{The first author acknowledges  the support of National
Natural Science Foundation of China (No.71561017), the Youth Academic Talent Plan of Lanzhou University of Finance and Economics. The second author acknowledges  Gansu Province Science Technology Plan (No.1606RJZA041). }}

\author{Jingjun Guo$^1$, Yanping Xiao$^{2},\thanks{Corresponding author:
hardxiao@126.com.}$\\
{\it\small 1.School of Statistics, Lanzhou University of Finance and Economics,} \\
{\it\small Lanzhou,730020,P.R.China }\\
{\it\small 2.School of Mathematics and Computer Science,} \\
{\it\small Northwest Minzu University,Lanzhou,730000,P.R.China} \\}
\date{}

\maketitle

\begin{abstract}
In this article,
existence of the  $k$-th order  derivatives
of  local time
$ \widehat{\alpha}^{(k)}(x,t)$ is considered for two
  d-dimensional fractional Ornstein-Uhlenbeck processes $X^{H_1}_t$ and
 $\widetilde{X}^{H_2}_s$ with Hurst parameters $H_1$ and
 $H_2$, respectively.  Moreover, H$\hat{o}$lder regularity condition of fractional Ornstein-Uhlenbeck process $X^{H}_t$ of local time
 $\tilde{\alpha}^{(k)}(x,t)$ is obtained   by some techniques using in Guo et al. (2017) and in Lou et al. (2017).



\bigskip

\end{abstract}

\section{Introduction}\label{sec-1}
In this article, we consider the following stochastic differential equation
\begin{equation}
\begin{array}{lll}
dX_t=-X_tdt+vdB_t^{H}, X_0=x,
\end{array}
\end{equation}
where $v$ is nonnegative coefficient, $B_t^H$ is a fractional Brownian motion. The solution $X_t^H$ of (1) is expressed as
\begin{equation}
\begin{aligned}
X_t=e^{-t}(x+v\int_0^te^sdB_s^{H}),
\end{aligned}
\end{equation}
which  is called a fractional Ornstein-Uhlenbeck  process (O-U, for short). $X_t^{H_1}$ and $\widetilde{X}_t^{H_2}$ denote fractional O-U processes for two independent fractional Bownian motions.

Local time of Gaussian process are important subjects  in probability theory
and their  derivatives have received  much attention  recently.
Guo et al. \cite{guo}, Jung et al. \cite{jung1} and \cite{jung2} respectively discussed exponent integral for derivative local time
of  fractional Brownian motion,
Tanaka formula and occupation-time formula.
On the other hand,  several  authors paid  attention to joint continuity and H$\hat{o}$lder
regularity  of  local time, see  e.g.,
Lou et al.  \cite{lou} in the case of  stochastic differential equation  and Ayache et al. \cite{ayache} in the case of fractional Brownian sheets.

Motivated by \cite{guo} and \cite{lou}, existence and H$\hat{o}$lder regularity of higher-order derivative    of local time for fractional O-U processes  are studied in this paper.
We are concerned with derivatives of  intersection local time of
$X^{H_1}$ and $\widetilde{X}^{H_2}$,   defined by
\begin{equation*}
\hat{\alpha}^{(k)}(x,t):= \frac{\partial ^k}{\partial x_1^{k_1}\cdots\partial x_d^{k_d}}
\int_0^t\int_0^t\delta(X^{H_1}_u-\widetilde{X}^{H_2}_s+x)duds,
\end{equation*}
where $k=(k_1, \cdots, k_d)$ is a multi-index with all
$k_i$ being nonnegative integers and $\delta$ is the Dirac delta
function of $d$-variable. In particular, local time of one fractional O-U process is given by
\begin{equation*}
\widetilde{\alpha}^{(k)}(x,t):= \frac{\partial ^k}{\partial x_1^{k_1}\cdots\partial x_d^{k_d}}
\int_0^t\delta(X^{H}_s+x)ds,
\end{equation*}
where $X^H$ is a fractional O-U process.
 Here $\delta^{(k)}(x)=\frac{\partial ^k}{\partial x_1^{k_1}\cdots\partial x_d^{k_d}}\delta(x)$ is k-th  order partial
derivative of the Dirac delta function.

Since the Dirac delta function $\delta$  is a generalized function,
we need  approximate the Dirac delta function $\delta$ by
\begin{equation}
f_{\varepsilon}(x,t):=\frac{1}{(2\pi\varepsilon)^{\frac{d}{2}}}
e^{-\frac{\mid x\mid^2}{2\varepsilon}}
=\frac{1}{(2\pi)^d}\int_{\mathbb{R}^d}e^{ipx}e^{-\frac{\varepsilon \mid p\mid^2}{2}}dp.
\end{equation}
Thus we approximate $\delta^{(k)}$ by
\begin{equation}
\begin{aligned}
f_{\varepsilon}^{(k)}(x,t):=\frac{\partial ^k}{\partial x_1^{k_1}\cdots\partial x_d^{k_d}} f_{\varepsilon}(x)
=\frac{i^k}{(2\pi)^d}\int_{\mathbb{R}^d}  p_1^{k_1}\cdots p_d^{k_d}
e^{ipx}e^{-\frac{\varepsilon
\mid p\mid^2}{2}}dp.
\end{aligned}
\end{equation}

\section{Existence of local time}\label{sec-2}
In this section, we discuss  existence of k-order derivative of two independent
 fractional O-U processes  by using similar method in Guo et al.(2017).

\begin{theorem}
 Let $B^{H_1}$ and $\widetilde{B}^{H_2}$   be two independent  d-dimensional fractional
 Brownian motions of Hurst parameters $H_1$ and $H_2$, respectively. Then local time  $\hat{\alpha}^{(k)}_{\varepsilon}(0,t)$ of fractional O-U
 processes $X^H_s$ and $\tilde{X}^H_s$ belongs to space $L^2(\Omega)$, when  $\varepsilon$ tend 0. Moreover, if
  its limit is denoted by $\hat{\alpha}^{(k)}(0,t)$, then $\hat{\alpha}^{(k)}(0,t)\in L^2(\Omega)$.
\end{theorem}

\noindent{\bf Proof}
Firstly, we claim that $\hat{\alpha}^{(k)}_{\varepsilon}(0,t)\in L^2(\Omega)$.
Indeed, denote $A_2=\{0<u,s<t\}^2$ and we have

\begin{eqnarray*}
\begin{aligned}
\mathbb{E}&\left[\left|\widehat{\alpha}^{(k)}_{\varepsilon}(0,t)\right|^2
\right]\\
&\le \frac{1}{(2\pi)^{2d}}\int_{A_2}\int_{\mathbb{R}^{2d}} \bigg| \mathbb{E}[\exp\{ip_1(X^{H_1}_{s_1}-\widetilde{X}^{H_2}_{u_1})+ip_2(X^{H_1}_{s_2}-\widetilde{X}^{H_2}_{u_2})\}]\bigg| \\
&~~\times\exp\{-\frac{\varepsilon}{2}\sum_{j=1}
^2\mid p_j\mid^2\}\prod_{j=1}^2
\mid p_j ^k\mid  dp dtds\\
&=\frac{1}{(2\pi)^{2d}}\int_{A_2}\int_{\mathbb{R}^{2d}}  \exp\{-\frac{1}{2}\mathbb{E}\big[\sum_{j=1}^2p_j(X^{H_1}_{s_j}-
\widetilde{X}^{H_2}_{t_j})
\big]^2\}
\\
&~~\times\exp\{-\frac{\varepsilon}{2}\sum_{j=1}
^n\mid p_j\mid^2\} \prod_{j=1}^2
\mid p_j ^k\mid  dpdtds\\
&\le \frac{1} {(2\pi)^{2d}}\int_{A_2}\int_{\mathbb{R}^{2d}}    \prod_{i=1}^d  \left(\prod_{j=1}^2
\mid p_{ij}^{k_i}\mid  \right)     \exp\{-\frac{1}{2}\mathbb{E}[p_{i1}X^{H_1,i}_{s_1}+p_{i2}X^{H_1,i}_{s_2}]^2\\
&~~-\frac{1}{2}\mathbb{E}[p_{i1}\widetilde{X}^{H_2,i}_{u_1}+p_{i2}\widetilde{X}^{H_2,i}_{u_2}]^2\}
dpduds \,.
\end{aligned}
\end{eqnarray*}
According to the fact in \cite{shen}, there are
\begin{eqnarray*}
Var\left(\xi(X_u^{H_1}-X_s^{H_1})\right)\geq C_1\xi^2(u-s)^{2H_1},
\end{eqnarray*}
and
\begin{eqnarray*}
Var\left(\eta(\widetilde{X}_u^{H_2}-\widetilde{X}_s^{H_2})\right)\geq C_2\eta^2(u-s)^{2H_2},
\end{eqnarray*}
where $C_1$ and $C_2$ are both some constants.

Similar method in \cite{guo}, we have
\begin{eqnarray*}
&&\mathbb{E} \left[\left| \widehat{\alpha}^{(k)}_{\varepsilon}(0,t)\right|^2
\right]
\le (2!)^2 C^2 \sum_{i,j=1}^2  \int_{A_2} (s_{i}-s_{ i-1 })^{- \rho H_1|k|} (u_{j}-u_{ j-1 })^{- (1-\rho)  H_2|k|}\\
  &&\qquad \qquad \qquad\qquad \cdot
  \left[ s_{1} (s_{2}-s_{1}) \right]^{-\gamma H_1d}\left[ u_{1} (u_{2}-u_{1}) \right]^{-(1-\gamma)H_2d }duds\,,
\end{eqnarray*}
where $A_2=\left\{ 0<s_1<s_2< t\right\}$ denotes the
simplex in $[0, t]^2$.   We choose $\rho=\gamma=\frac{H_2}{H_1+H_2}$  to obtain

\begin{eqnarray*}
&&\mathbb{E} \left[\left| \widehat{\alpha}^{(k)}_{\varepsilon}(0,t)\right|^2
\right]
  \le C_0(2!)^{2-2\kappa}t^{4\kappa},
\end{eqnarray*}
where $C_0$ is a constant independent of $t$ and $\kappa$ is some constant. Thus, the left of (5)
is finite if $\frac{H_1H_2}{H_1+H_2}(|\kappa|+d)\leq1$.

Secondly, we claim that the sequence $\{\widehat{\alpha}^{(k)}_{\varepsilon}(0,t),\varepsilon>0\}$ is a Cauchy sequence in $L^2$. In fact, for any $\varepsilon,\theta>0$, there is
\begin{eqnarray*}
\begin{aligned}
\mathbb{E}&\left[\left|\widehat{\alpha}^{(k)}_{\varepsilon}(0,t)-
\widehat{\alpha}^{(k)}_{\theta}(0,t)\right|^2
\right]\\
&\le \frac{1}{(2\pi)^{2d}}\int_{A_2}\int_{\mathbb{R}^{2d}} \bigg| \mathbb{E}[\exp\{ip_1(X^{H_1}_{s_1}-\widetilde{X}^{H_2}_{u_1})+ip_2(X^{H_1}_{s_2}-
\widetilde{X}^{H_2}_{u_2})\}]\bigg|
\\
&~~\times\bigg|\exp\{-\frac{\varepsilon}{2}\sum_{j=1}
^2\mid p_j\mid^2\}-\exp\{-\frac{\theta}{2}\sum_{j=1}
^2\mid p_j\mid^2\}\bigg|\prod_{j=1}^2
\mid p_j ^k\mid  dp dtds\\
&=\frac{1}{(2\pi)^{2d}}\int_{A_2}\int_{\mathbb{R}^{2d}}  \exp\{-\frac{1}{2}\mathbb{E}\big[\sum_{j=1}^2p_j(X^{H_1}_{s_j}-
\widetilde{X}^{H_2}_{t_j})
\big]^2\}
\\
&~~\times\left[1-\exp\{-\frac{|\varepsilon-\theta|}{2}\sum_{j=1}
^n\mid p_j\mid^2\}\right]\prod_{j=1}^2
\mid p_j ^k\mid  dpdtds\\
&\le \frac{\sup_{p_j\in R}\left\{1-\exp\{-\frac{|\varepsilon-\theta|}{2}\sum_{j=1}
^n\mid p_j\mid^2\}\right\}} {(2\pi)^{2d}}\int_{A_2}\int_{\mathbb{R}^{2d}}    \prod_{i=1}^d  \left(\prod_{j=1}^2
\mid p_{ij}^{k_i}\mid  \right)\\
&~~\times\exp\{-\frac{1}{2}\mathbb{E}[p_{i1}X^{H_1,i}_{s_1}+p_{i2}X^{H_1,i}_{s_2}]^2-\frac{1}{2}\mathbb{E}[p_{i1}\widetilde{X}^{H_2,i}_{u_1}+p_{i2}\widetilde{X}^{H_2,i}_{u_2}]^2\}
dpduds \,.
\end{aligned}
\end{eqnarray*}
According to dominated convergence theorem, we obtain
\begin{eqnarray*}
&&\mathbb{E} \left[\left| \widehat{\alpha}^{(k)}_{\varepsilon}(0,t)-\widehat{\alpha}^{(k)}_{\theta}(0,t)\right|^2
\right]\rightarrow 0,
\end{eqnarray*}
as $\varepsilon\rightarrow 0$ and $\theta\rightarrow 0.$ Hence, $\{\widehat{\alpha}^{(k)}_{\varepsilon}(0,t),\varepsilon>0\}$ is a Cauchy sequence in $L^2(\Omega)$, which means that $\widehat{\alpha}^{(k)}(0,t)$ belongs to space  $L^2(\Omega)$.
\hfill $\Box$

\noindent{\bf Remark} Comparing this condition in Theorem 1 with Guo et al. (2017), we
find that there is the same condition, which is $\frac{H_1H_2}{H_1+H_2}(\mid \kappa+d\mid)\leq1.$
That is said that there have the same condition to exist in space $L^2(\Omega)$.

\section{The H$\hat{o}$lder regularity}\label{sec-3}
In this section, we discuss  condition of H$\hat{o}$lder regularity of higher-order derivative of fractional O-U process $X^{H}_t$, which satisfies a stochastic differential equation driven by fractional Brownian motion $B^H_t$. Firstly, we give definition of  H$\hat{o}$lder exponent of local time $\widetilde{\alpha}^{(k)}(t,x)$. Secondly, condition of  H$\hat{o}$lder regularity is  obtained through some necessary lemmas.

\begin{definition}
The pathwise H$\hat{o}$lder exponent of local time $\widetilde{\alpha}^{(k)}(x,t)$ is defined by
\begin{eqnarray*}
\alpha(t)=\sup\{\alpha>0,\limsup_{h\rightarrow 0}\sup_{x\in\mathbb{R}^d}\frac{\widetilde{\alpha}^{(k)}(x,t+h)-\widetilde{\alpha}^{(k)}(x,t)}
{h^{\alpha}}=0\}.
\end{eqnarray*}
\end{definition}
To state our main result in this section, we need some necessary lemmas.
\begin{lemma}
 Let $n\geq1$ and $k\geq2$  be arbitrary. If dimensional d and Hurst parameter H satisfy  $H(\mid k\mid+d)\leq1$, then
 \begin{eqnarray*}
\begin{aligned}
\mathbb{E}\left[\left|\widetilde{\alpha}^{(k)}_{\varepsilon}(x,t+h)-
\widetilde{\alpha}^{(k)}_{\varepsilon}(x,t)\right|^n
\right]\le C_3 h^{n-nHd-H\mid k\mid},
\end{aligned}
\end{eqnarray*}
where $C_3=\frac{nn!C^n}{\Gamma(n-nHd-H\mid k\mid+1)}.$
\end{lemma}

\noindent{\bf Proof}
The similar techniques in Guo et al.(2017) is used to prove the lemma. Fix an integer $n\geq 1$. Denote  $B_n=\{t<u,s <t+h\}^n$.
We have
\begin{eqnarray*}
\begin{aligned}
\mathbb{E}&\left[\left|\widetilde{\alpha}^{(k)}_{\varepsilon}(x,t+h)-
\widetilde{\alpha}^{(k)}_{\varepsilon}(x,t)\right|^n
\right]\\
&\le \frac{1} {(2\pi)^{nd}}\int_{B_n}\int_{\mathbb{R}^{nd}}    \prod_{i=1}^d  \left(\prod_{j=1}^n
\mid p_{ij}^{k_i}\mid  \right)     \exp\{-\frac{1}{2}\mathbb{E}[p_{i1}X^{i}_{s_1}+\cdots+p_{in}X^{i}_{s_n}]^2\}dpds \,.
\end{aligned}
\end{eqnarray*}
  The expectation in the
above exponent  can be computed by
\begin{eqnarray*}
&&\mathbb{E}[p_{i1}X^{i}_{s_1}+\cdots+p_{in}X^{i}_{s_n}]^2=
(p_{i1}, \cdots, p_{in}) Q(p_{i1}, \cdots, p_{in})^T,
\end{eqnarray*}
where
\[
Q=\EE \left(X^{i}_j X^{i}_k\right)_{1\le j,k\le n}
\quad\]
 denotes  covariance matrix  of n-dimensional random vector
$(X^{i}_{s_1},...,X^{i}_{s_n})$.  Thus we have
\begin{eqnarray*}
\mathbb{E}\left[
\left|\widetilde{\alpha}^{(k)}_{\varepsilon}(x,t+h)-
\widetilde{\alpha}^{(k)}_{\varepsilon}(x,t)\right|^n
\right]\le \frac{1}{(2\pi)^{nd}}  \int_{B_n}  \prod_{i=1}^d  I_i(s)ds\,,
\label{e.boundedbyI}
\end{eqnarray*}
where
\[
I_i(s):= \int_{\mathbb{R}^{n}}
\mid x^{k_i} \mid \exp\{-\frac{1}{2}x^TQx\} dx  \,.
\]
Here we recall $x=(x_1, \cdots, x_n)$  and
$x^k_i=x_1^{k_i} \cdots x_n^{k_i}$.

  Making substitution
$\xi=\sqrt{Q}x$, then
\begin{eqnarray*}
I_i(s)=\int_{\mathbb{R}^{n }}\prod_{j=1}^n
\mid (Q^{-\frac{1}{2}}\xi)_j\mid^{k_i} \exp\{-\frac{1}{2}\mid \xi\mid^2\}\mathrm{det}(Q)
^{-\frac{1}{2}}d\xi.\\
\end{eqnarray*}
To obtain a nice bound for the above
integral, let us first diagonalize   $Q$:
\begin{equation*}
Q=P\Lambda P^{-1} \,,
\end{equation*}
where $\Lambda=$diag$\{\lambda_1 ,...,\lambda_n \}$ is a strictly
positive diagonal matrix with $\la_1\le \la_2\le\cdots\le \la_d$
and    $P=(q_{ij})_{1\le i,j\le d} $ is  an orthogonal matrix.  Hence, we have
$\det(Q)=\la_1\cdots\la_d$.
Denote
\begin{equation*}
\eta=
\begin{pmatrix}
\eta_1, \eta_2, \cdots,
\eta_n
\end{pmatrix}^T=P^{-1}\xi.
\end{equation*}
Hence,
\begin{eqnarray*}
P^{-\frac{1}{2}}\xi
&=&P\Lambda^{-1/2}P^{-1} \xi =
P\Lambda^{-1/2}\eta\\
&=&P\begin{pmatrix}
\lambda_1^{-\frac{1}{2}}\eta_1\\\lambda_2^{-\frac{1}{2}}\eta_2\\\vdots\\\lambda_n^{-\frac{1}{2}}
\eta_n
\end{pmatrix}
=
\begin{pmatrix}
q_{1,1}&q_{1,2}&\cdots&q _{1,n} \\
q_{2,1}&q_{2,2}&\cdots&q_{2,n}\\
\vdots&\vdots&\cdots&\vdots\\
q_{n,1}&q _{n,2}&\cdots&q_{n,n}\\
\end{pmatrix}
\begin{pmatrix}
\lambda_1^{-\frac{1}{2}}\eta_1\\\lambda_2^{-\frac{1}{2}}\eta_2\\\vdots\\\lambda_n^{-\frac{1}{2}}
\eta_n
\end{pmatrix}.
\end{eqnarray*}

Therefore, we have
\begin{eqnarray*}
\mid(Q^{-\frac{1}{2}}\xi)_j\mid
&=&\mid\sum_{k=1}^nq_{jk}\lambda_{k}^{-\frac{1}{2}}\eta_{k}\mid
\leq
\lambda_1
^{-\frac{1}{2}}  \sum_{k=1}^n\mid q_{jk} \eta_{k}\mid \\
&\leq & \lambda_1
^{-\frac{1}{2}}
\left(\sum_{k=1}^nq_{jk}^2\right)^{\frac{1}{2}}
\left(\sum_{k=1}^n\eta_{k}^2\right)^{\frac{1}{2}}
 \leq  \lambda_1
^{-\frac{1}{2}} \mid\eta\mid_2
 =\lambda_1
^{-\frac{1}{2}}\mid\xi\mid_2.\\
\end{eqnarray*}
 Since  $Q$ is positive definite, we see that
 \[
 \la_1\ge \la_1(Q),
 \]
 where $  \la_1(Q)$ is the smallest eigenvalue of
 $Q$. Hence, we have
\begin{eqnarray*}
I_i(s)= \mathrm{det}(Q)
^{-\frac{1}{2}} \la_1(Q)^{-\frac12 k_i }
\int_{\mathbb{R}^{n }}
\mid\xi\mid_2 ^{k_i}  \exp\{-\frac{1}{2}\mid \xi\mid^2\} d\xi.
\end{eqnarray*}

Now we are going to find a lower bound for $\la_1(Q_1)$. By the well-known result
\begin{eqnarray*}
\lambda_1(Q)&\geq & K\min\left[ s_{1}^{2H}, (s_{2}-s_{1})^{2H}, \cdots ,(s_{n}-s_{n-1})^{2H} \right],
\end{eqnarray*}
we have
\begin{eqnarray*}
&\int_{\mathbb{R}^n}\prod_{j=1}^n\mid(Q^{-\frac{1}{2}\xi})_j\mid^{k_i}
\exp\{-\frac{1}{2}\mid\xi\mid^2\}det(Q)^{-\frac{1}{2}}d\xi\\
&\le  C^n\min_{j=1,...,n}(s_j-s_{j-1})^{-Hk_i}\left[ s_{1}(s_{2}-s_{1})\cdots (s_{n}-s_{n-1})^{-H} \right].\\
\end{eqnarray*}

Therefore,
\begin{eqnarray*}
&&\mathbb{E} \left[\left| \widetilde{\alpha}^{(k)}_{\varepsilon}(x,t+h)-
\widetilde{\alpha}^{(k)}_{\varepsilon}(x,t)\right|^n
\right]\\
&&\le n! C^n  \int_{B_n}\min_{j=1,...,n}(s_{j}-s_{ j-1 })^{-  H|k|}
  \left[ s_{1} (s_{2}-s_{1}) \cdots (s_{n}-s_{n-1}) \right]^{-Hd}ds\\
&&\le n! C^n \sum_{i=1}^n  \int_{B_n} (s_{i}-s_{ i-1 })^{- H|k|}\left[ s_{1} (s_{2}-s_{1}) \cdots (s_{n}-s_{n-1}) \right]^{-Hd}ds\,.
\end{eqnarray*}
If $H(|k|+d)\leq1$, then
\begin{eqnarray*}
&&\int_{A_n}(s_{i}-s_{ i-1 })^{-  H|k|}
  \left[ s_{1} (s_{2}-s_{1}) \cdots (s_{n}-s_{n-1}) \right]^{-Hd}ds\\
&&\le \frac{C^n h^{n-nHd-H|k|}}{\Gamma(n-nHd-H|k|+1)}.
\end{eqnarray*}
Thus,
\begin{eqnarray*}
\begin{aligned}
\mathbb{E}\left[\left|\widetilde{\alpha}^{(k)}_{\varepsilon}(x,t+h)-
\widetilde{\alpha}^{(k)}_{\varepsilon}(x,t)\right|^n
\right]\le C_3 h^{n-nHd-H\mid k\mid}.
\end{aligned}
\end{eqnarray*}
\hfill $\Box$

\noindent{\bf Remark} By Lemma 1, there is the following inequality
\begin{eqnarray*}
\begin{aligned}
\mathbb{E}\left[\left|\frac{\widetilde{\alpha}^{(k)}_{\varepsilon}(x,t+h)-
\widetilde{\alpha}^{(k)}_{\varepsilon}(x,t)
}{h^{1-Hd}}\right|^n\right]\le C_4(n!)^{2-Hd},
\end{aligned}
\end{eqnarray*}
where $C_4$ is a constant.

\begin{lemma}
 Let $n\geq1$ and $k\geq2$  be arbitrary. If dimensional $d\geq1$, any $\delta$ and Hurst parameter H satisfying  $H(\mid k+\delta\mid+d)\leq1$, then
 \begin{eqnarray*}
\begin{aligned}
\mathbb{E}\left[\left|\widetilde{\alpha}^{(k)}(x,t)-
\widetilde{\alpha}^{(k)}(y,t)\right|^n
\right]\le C_5\mid x-y\mid^{\delta},
\end{aligned}
\end{eqnarray*}
where $\mid x-y\mid^{\delta}\equiv\sum_{j=1}^n|x_j-y_j|^{\delta_j}$ and $C_5=\frac{C^nn!}{(2\pi)^{nd}}\cdot\frac{t^{n-nHd-H|k+\delta|}}{\Gamma(n-nHd-H|k+\delta|+1)}$.
\end{lemma}

\noindent{\bf Proof} We only verify that
\begin{eqnarray*}
\begin{aligned}
\mathbb{E}\left[\left|\widetilde{\alpha}^{(k)}_{\varepsilon}(x,t)-
\widetilde{\alpha}^{(k)}_{\varepsilon}(y,t)\right|^n
\right]\le C_5\mid x-y\mid^{\delta}.
\end{aligned}
\end{eqnarray*}
Indeed according to the definition of $\widetilde{\alpha}^{(k)}_{\varepsilon}(x,t)$, we have
\begin{eqnarray*}
\begin{aligned}
\mathbb{E}&\left[\left|\widetilde{\alpha}^{(k)}_{\varepsilon}(x,t)-
\widetilde{\alpha}^{(k)}_{\varepsilon}(y,t)\right|^n
\right]\\
&=E\left[\left|\frac{i^kp^k}{(2\pi)^d}\int_{\mathbb{R}^d}\int_0^t\left(e^{ip(X_s^H-x)}-e^{ip(X_s^H-y)} \right)e^{-\frac{\varepsilon}{2}|p|^2}dpds\right|^n\right]\\
&=E\left[\left|\frac{i^kp^k}{(2\pi)^d}\int_{\mathbb{R}^d}\int_0^te^{ipX_s^H}
\left(e^{-ipx}-
e^{-ipy}\right)e^{-\frac{\varepsilon}{2}|p|^2}dpds\right|^n\right].
\end{aligned}
\end{eqnarray*}
Using the well-known  inequality
\begin{eqnarray*}
\begin{aligned}
\left|e^{-ipa}-e^{-ipb}\right|\leq C_6|p|^{\delta}|a-b|^{\delta},\\
\end{aligned}
\end{eqnarray*}
where $\delta$ and $C_6$  are both some constants, there is

\begin{eqnarray*}
\begin{aligned}
\mathbb{E}&\left[\left|\widetilde{\alpha}^{(k)}_{\varepsilon}(x,t)-
\widetilde{\alpha}^{(k)}_{\varepsilon}(y,t)\right|^n
\right]\\
&\le \frac{1} {(2\pi)^{nd}}\int_{A_n}\int_{\mathbb{R}^{nd}} \left|E\left[\prod_{j=1}^n
\mid p_{j}^{k}\mid e^{ip_jX_{s_j}^{H}}\left(e^{-ip_jx_j}-e^{-ip_jy_j}\right)
\right]\right|dpds\\
&\le\frac{C_6\prod_{j=1}^d\left|x_j-y_j\right|^{\delta_j}} {(2\pi)^{nd}}\int_{A_n}\int_{\mathbb{R}^{nd}} \prod_{j=1}^n
\mid p_{j}\mid^{\delta_j+k} \left|E\left[\exp\{ip_1X_{s_1}^H+...+ip_nX_{s_n}^H\}\right]\right|\\
&\qquad\cdot\exp\{-\frac{\varepsilon}{2}\sum_{j=1}^d
|p_j|^2\}dpds.\\
\end{aligned}
\end{eqnarray*}
  The expectations in the
above exponent  can be computed by
\begin{eqnarray*}
\mathbb{E}[p_{i1}X^{H,i}_{s_1}+\cdots+p_{in}X^{H,i}_{s_n}]^2=
(p_{i1}, \cdots, p_{in}) Q(p_{i1}, \cdots, p_{in})^T,
\end{eqnarray*}
where
\[
Q=\EE \left(X^{H, i}_j X^{H, i}_k\right)_{1\le j,k\le n}
,\]
 denotes covariance matrix  of n-dimensional random vector
$(X^{H,i}_{s_1},...,X^{H,i}_{s_n})$.  Thus, we obtain
\begin{eqnarray}
\mathbb{E}\left[
\left|\widetilde{\alpha}^{(k)}_{\varepsilon}(x,t)-\widetilde{\alpha}^{(k)}_{\varepsilon}(y,t)\right|^n
\right]\le \frac{C_6}{(2\pi)^{nd}} \prod_{i=1}^d|x_i-y_i|^{\delta_i}\int_{A_n}  \prod_{i=1}^d  I_i(s)ds\,,
\label{e.boundedbyI}
\end{eqnarray}
where
\[
I_i(s):= \int_{\mathbb{R}^{n}}
\mid p_i^{k} \mid^{\delta_i} \exp\{-\frac{1}{2}P^TQP\} dP.
\]
Here we recall $P=(p_1, \cdots, p_n)$  and
$p^k_i=p_1^{k_i} \cdots p_n^{k_i}$.

Making substitution
$\xi=\sqrt{Q}p$, then
\begin{eqnarray*}
I_i(s)=\int_{\mathbb{R}^{n }}\prod_{j=1}^n
\mid (Q^{-\frac{1}{2}}\xi)_j\mid^{k_i+\delta_i} \exp\{-\frac{1}{2}\mid \xi\mid^2\}\mathrm{det}(Q)
^{-\frac{1}{2}}d\xi.\\
\end{eqnarray*}
Let us first diagonalize   $Q$:
\begin{equation*}
Q =R\Lambda R^{-1} \,,
\end{equation*}
where $\Lambda=$diag$\{\lambda_1 ,...,\lambda_n \}$ is a strictly
positive diagonal matrix with $\la_1\le \la_2\le\cdots\le \la_d$
and    $R =(r_{ij})_{1\le i,j\le d} $ is  an orthogonal matrix.  Hence, we have
$\det(Q)=\la_1\cdots\la_d$.
Denote
\begin{equation*}
\eta=
\begin{pmatrix}
\eta_1, \eta_2, \cdots,
\eta_n
\end{pmatrix}^T=R^{-1}\xi.
\end{equation*}
Hence,
\begin{eqnarray*}
R^{-\frac{1}{2}}\xi
&=&R\Lambda^{-1/2}R^{-1} \xi =
R\Lambda^{-1/2}\eta\\
&=&R\begin{pmatrix}
\lambda_1^{-\frac{1}{2}}\eta_1\\\lambda_2^{-\frac{1}{2}}\eta_2\\\vdots\\\lambda_n^{-\frac{1}{2}}
\eta_n
\end{pmatrix}
=
\begin{pmatrix}
r_{1,1}&r_{1,2}&\cdots&r _{1,n} \\
r_{2,1}&r_{2,2}&\cdots&r_{2,n}\\
\vdots&\vdots&\cdots&\vdots\\
r_{n,1}&r _{n,2}&\cdots&r_{n,n}\\
\end{pmatrix}
\begin{pmatrix}
\lambda_1^{-\frac{1}{2}}\eta_1\\\lambda_2^{-\frac{1}{2}}\eta_2\\\vdots\\\lambda_n^{-\frac{1}{2}}
\eta_n.
\end{pmatrix}.
\end{eqnarray*}

Hence, we have
\begin{eqnarray*}
\mid(Q^{-\frac{1}{2}}\xi)_j\mid
 \leq  \lambda_1
^{-\frac{1}{2}} \mid\eta\mid_2
 =\lambda_1
^{-\frac{1}{2}}\mid\xi\mid_2.\\
\end{eqnarray*}
 Since  $Q$ is positive definite, we see that
 \[
 \la_1\ge \la_1(Q)
 \]
 where $  \la_1(Q)$ is the smallest eigenvalue of
 $Q$.
 This implies
 \begin{equation*}
\mid(Q^{-\frac{1}{2}}\xi)_j\mid
\le \la_1(Q)^{-\frac12}
\mid\xi\mid_2.\label{e.Bxi_bound}
\end{equation*}
Consequently, we have
\begin{eqnarray*}
I_i(s)&=& \mathrm{det}(Q)
^{-\frac{1}{2}} \la_1(Q)^{-\frac12(k_i+\delta_i)}\int_{\mathbb{R}^{n }}
\mid\xi\mid_2 ^{k_i+\delta_i}  \exp\{-\frac{1}{2}\mid \xi\mid^2\} d\xi.
\label{e.I_first_bound}
\end{eqnarray*}

Next we need find a lower bound for $\la_1(Q)$ as following
\begin{eqnarray*}
\la_1(Q)
&\ge&  K \min\{s_{1}^{2H},(s_{2}-s_{1})^{2H},...,
(s_{n}-s_{n-1})^{2H}\}\,. \label{e.la_1_q2}
\end{eqnarray*}
Using the following fact
\begin{equation*}
\begin{aligned}
\mathrm{det}(Q)
\geq& C_7^ns_{1}^{2H}(s_{2}-s_{1})^{2H}\cdots (s_{n}-s_{n-1})^{2H},\\
\end{aligned}
\end{equation*}
 we  have
\begin{eqnarray*}
&&\mathbb{E} \left[\left| \widetilde{\alpha}^{(k)}_{\varepsilon}(x,t)-\widetilde{\alpha}^{(k)}_{\varepsilon}(y,t)\right|^n
\right] \le \frac{n!C_7^n}{(2\pi)^{nd}} \prod_{j=1}^n(x_j-y_j)^{\delta_j} \\
&&\qquad \times \sum_{j=1}^n\int_{A_n}(s_j-s_{j-1})^{-H|k+\delta|}\left[s_1(s_2-s_1)...(s_n-s_{n-1})\right]^{-Hd}ds.\\
\end{eqnarray*}

If $H(|k+\delta|+d)\leq1$, then
\begin{eqnarray*}
&&\int_{A_n}(s_j-s_{j-1})^{-H|k+\delta|}\left[s_1(s_2-s_1)...(s_n-s_{n-1})\right]^{-Hd}ds,\\
&&\leq \frac{n! C_7^n }{(2\pi)^{nd}}\prod_{j=1}^n(x_j-y_j)^{\delta_j} \frac{t^{n-nHd-H|k+\delta|}}{\Gamma(n-nHd-H|k+\delta|+1)}\\
&&\equiv C_8\prod_{j=1}^n(x_j-y_j)^{\delta_j}, \\
\end{eqnarray*}
where $C_8=\frac{n! C_7^n }{(2\pi)^{nd}}\times\frac{t^{n-nHd-H|k+\delta|}}{\Gamma(n-nHd-H|k+\delta|+1)}$.
So this lemma is proved.
\hfill $\Box$

\begin{theorem}
The pathwise H$\hat{o}$lder exponent of local time $\widetilde{\alpha}^{(k)}(x,t)$ is given by
\begin{eqnarray*}
\alpha(t)=n-nHd-H|k|.
\end{eqnarray*}
\end{theorem}

\noindent{\bf Proof}
According to Lemma 1,2,
 there is
\begin{eqnarray*}
\alpha(H)\geq n-nHd-H|k|.
\end{eqnarray*}
Next to establish this theorem, we only prove that
\begin{eqnarray*}
\alpha(H)\leq n-nHd-H|k|.
\end{eqnarray*}
Indeed, fix any $t>0,$ since the property of local time, we have
\begin{equation*}
\begin{array}{lll}
h&=\int_{\mathbb{R}^d}\left(\widetilde{\alpha}^{(k)}(x,t+h)-\widetilde{\alpha}^{(k)}(x,t)\right)dx\\
&\leq\sup_{x\in\mathbb{R}^d}\left\{\widetilde{\alpha}^{(k)}(x,t+h)-\widetilde{\alpha}^{(k)}(x,t)
\right\}
\sup_{s,u\in[t,t+h]}\mid X^H_u-X^H_s\mid\\
&\leq2\sup_{x\in\mathbb{R}^d}\left\{\widetilde{\alpha}^{(k)}(x,t+h)-\widetilde{\alpha}^{(k)}(x,t)
\right\}
\sup_{s\in[t,t+h]}\mid X^H_t-X^H_s\mid.
\end{array}
\end{equation*}

By Yan et al. (2008), fractional O-U process $X^H_t$ can be written by
\begin{eqnarray*}
X^H_t=v\int_0^tF(t,u)dB_u,
\end{eqnarray*}
where
\begin{eqnarray*}
F(t,u)=(H-\frac{1}{2})K_He^{-t}u^{\frac{1}{2}-H}\int_u^ts^{H-\frac{1}{2}}(s-u)
^{H-\frac{3}{2}}e^sds,
\end{eqnarray*}
if $\frac{1}{2}<H<1$,
\begin{equation*}
\begin{array}{lll}
F(t,u)&=&K_Hu^{\frac{1}{2}-H}(-e^{-t}\int_u^t(s-u)
^{H-\frac{1}{2}}s^{H-\frac{1}{2}}e^sds\\
&&+t^{H-\frac{1}{2}}(t-u)^{H-\frac{1}{2}}+\frac{2}{1-2H}e^{-t}
\int_u^t(s-u)^{H-\frac{1}{2}}s^{H-\frac{3}{2}}e^sds,
\end{array}
\end{equation*}
if $0<H<\frac{1}{2}$, with $K_H=\frac{2H\Gamma(\frac{3}{2}-H)}{\Gamma(H+\frac{1}{2})\Gamma(2-2H)^{\frac{1}{2}}}$,
 then
\begin{equation*}
\begin{array}{lll}
\mid X^H_t-X^H_s\mid&=\mid v\int_0^tF(t,u)dB_u-v\int_0^sF(s,u)dB_u\mid\\
&\leq\mid v\mid\sup\mid F(t,u)\mid\cdot\mid B_t-B_s\mid\\
&\leq\mid v\mid\sup\mid F(t,u)\mid\cdot\max\{\mid B_{t-s}\mid\}\cdot\mid t-s\mid.
\end{array}
\end{equation*}

Therefore,
\begin{eqnarray*}
&&h\leq2\mid v\mid\sup_{x\in\mathbb{R}^d}\left(\widetilde{\alpha}^{(k)}(x,t+h)-\widetilde{\alpha}^{(k)}(x,t)\right)
\sup_{s\in[t,t+h]}\mid F(t,u)\mid\max\{\mid B_{t-s}\mid\}\mid t-s\mid.
\end{eqnarray*}
Using Lemma 1, we get
\begin{eqnarray*}
\frac{\sup_{x\in\mathbb{R}^d}\left\{\widetilde{\alpha}^{(k)}(x,t+h)-\widetilde{\alpha}^{(k)}(x,t)\right\}}
{h^{n-nHd-H\mid k\mid}}\geq C_9,
\end{eqnarray*}
where $C_9=\frac{1}{2|v|\sup|F(t,u)|\max\{|B_{t-s}|\}}.$

\hfill $\Box$

\end{document}